\documentclass{amsart}

 \theoremstyle{definition}
 
 \theoremstyle{remark}
 
 \numberwithin{equation}{subsection}
\begin{document}

\title[2-LOCAL DERIVATIONS ON VON NEUMANN ALGEBRAS ]
 {2-LOCAL DERIVATIONS ON VON NEUMANN ALGEBRAS OF TYPE I}

\author{ Shavkat Ayupov}

\address{Institute of Mathematics and Information Technologies,
Tashkent, Uzbekistan and the Abdus Salam International Centre for Theoretical Physics (ICTP) Trieste, Italy}

\email{sh$_-$ayupov@mail.ru}

\author{ Farkhad Arzikulov}

\address{Institute of Mathematics and Information Technologies,
Tashkent, Andizhan State University, Andizhan, Uzbekistan}

\email{arzikulovfn@rambler.ru}

%\thanks{This work was completed with the support of an Izaak
% Walton Killam Memorial Scholarship.}

%\thanks{The author was also supported in part by the Research
% Council of Slovenia.}

%\subjclass{Primary 46L57; Secondary 46L40}


\date{February 21, 2012.}

\dedicatory{}

%\commby{Daniel J. Rudolph}

%%% ----------------------------------------------------------------------

\begin{abstract}
In the present paper we prove that every 2-local derivation on a
von Neumann algebra of type I is a derivation.
\end{abstract}

%%% ----------------------------------------------------------------------
\maketitle
%%% ----------------------------------------------------------------------

\section*{Introduction}

The present paper is devoted to 2-local derivations on von Neumann
algebras. Recall that a 2-local derivation is defined as follows:
given an algebra $A$, a map $\bigtriangleup : A \to A$ (not linear
in general) is called a 2-local derivation if for every $x$, $y\in
A$, there exists a derivation $D_{x,y} : A\to A$ such that
$\bigtriangleup(x)=D_{x,y}(x)$ and $\bigtriangleup(y)=D_{x,y}(y)$.

In 1997, P. \v{S}emrl \cite{S} introduced the notion of 2-local derivations and
described 2-local derivations on the algebra $B(H)$ of all
bounded linear operators on the infinite-dimensional separable
Hilbert space H. A similar description for the finite-dimensional case appeared
later in \cite{KK}. In the paper \cite{LW} 2-local derivations have
been described on matrix algebras over finite-dimensional division rings.

In \cite{AK} the authors suggested a new technique and have generalized the above
mentioned results of \cite{S} and \cite{KK} for arbitrary Hilbert spaces. Namely they considered
2-local derivations on the algebra $B(H)$ of all linear
bounded operators on an arbitrary (no separability is assumed) Hilbert space $H$ and
proved that every 2-local derivation on $B(H)$ is a derivation.

In the present paper we also suggest another technique and
generalize the above mentioned results of \cite{S}, \cite{KK} and
\cite{AK} for arbitrary von Neumann algebras of type I. Namely, we
prove that every 2-local derivation on a von Neumann algebra of
type I is a derivation.

The authors want to thank K.K.Kudaybergenov for many stimulating conversations
on the subject.

\section{Preliminaries}

Let $M$ be a von Neumann algebra.

{\it Definition.} A linear map $D : M\to M$ is called a
derivation, if $D(xy)=D(x)y+xD(y)$ for any two elements $x$, $y\in
M$.

A map $\Delta : M\to M$ is called a 2-local derivation, if
for any two elements $x$, $y\in M$ there exists a derivation
$D_{x,y}:M\to M$ such that $\Delta (x)=D_{x,y}(x)$, $\Delta
(y)=D_{x,y}(y)$.

It is known that any derivation $D$ on a von Neumann algebra $M$ is an inner
derivation, that is there exists an element $a\in M$ such that
$$
D(x)=ax-xa, x\in M.
$$
Therefore for a von Neumann algebra $M$ the above definition is equivalent to the
following one: A map $\Delta : M\to M$ is called a
2-local derivation, if for any two elements $x$, $y\in M$
there exists an element $a\in M$ such that $\Delta (x)=ax-xa$,
$\Delta (y)=ay-ya$.

Further we will use the latter definition.

Let $n$ be an arbitrary infinite cardinal number,
$\Xi$ be a set of indexes of the cardinality $n$. Let $\{e_{ij}\}$
be a set of matrix units such that $e_{ij}$ is a $n\times
n$-dimensional matrix, i.e.
$e_{ij}=(a_{\alpha\beta})_{\alpha\beta\in\Xi}$, the $(i,j)$-th
component of which is $1$, i.e. $a_{ij}=1$, and the rest
components are zeros. Let $\{m_\xi\}_{\xi\in \Xi}$  be a set of
$n\times n$-dimensional matrixes. By $\sum_{\xi\in \Xi} m_\xi$ we
denote the matrix whose components are sums of the corresponding
components of matrixes of the set $\{m_\xi \}_{\xi\in \Xi}$. Let
$$
M_n({\bf C})=\{\{\lambda_{ij}e_{ij}\}: \,for\,\, all\,\,
indexes\,\, i,\,j \,\lambda_{ij}\in {\bf C},
$$
$$
and\,\, there\,\, exists\,\, such\,\, number\,\, K\in {\bf
R},\,\,that \,\, for \,\, all\,\, n\in N
$$
$$
and\,\, \{e_{kl}\}_{kl=1}^n\subseteq \{e_{ij}\} \Vert\sum_{kl=1}^n
\lambda_{kl}e_{kl}\Vert \le K\},
$$
where $\Vert \,\, \Vert$ is a norm of a matrix. It is easy to see
that $M_n({\bf C})$ is a vector space.

The associative multiplication of elements in $M_n({\bf C})$ can
be defined as follows: if $x=\sum_{ij\in \Xi}\lambda_{ij}e_{ij}$,
$ y=\sum_{ij\in \Xi}\mu_{ij}e_{ij}$ are elements of $M_n({\bf C})$
then $xy=\sum_{ij\in \Xi} \sum_{\xi\in \Xi} \lambda_{i\xi}\mu_{\xi
j}e_{ij}$. With this operation $M_n({\bf C})$ becomes an associative
algebra and $M_n({\bf C})=B(l_2(\Xi))$, where $l_2(\Xi)$ is a Hilbert space over ${\bf C}$ with elements
$\{x_i\}_{i\in \Xi}$, $x_i\in \bf C$ for all $i\in \Xi$, $B(l_2(\Xi))$ is the
associative algebra of all bounded linear operators on the Hilbert
space $l_2(\Xi)$. Then $M_n({\bf C})$ is a von Neumann algebra of
infinite $n\times n$-dimensional matrices over ${\bf C}$.

Similarly, if we take the algebra $B(H)$ of all bounded linear
operators on an arbitrary Hilbert space $H$ and if $\{q_i\}$ is an
arbitrary maximal orthogonal set of minimal projections of the
algebra $B(H)$, then $B(H)=\sum_{ij}^\oplus q_i B(H)q_j$ (see
\cite{AFN}).

Let $X$ be a hyperstonean compact, and let $C(X)$ denote the
commutative algebra of all complex-valued continuous functions on
the compact $X$ and
$$
\mathcal{M}=\{\{\lambda_{ij}(x)e_{ij}\}_{ij\in\Xi}: (\forall ij\,\,\,
\lambda_{ij}(x)\in C(X))
$$
$$
(\exists K\in R)(\forall m\in N)(\forall
\{e_{kl}\}_{kl=1}^m\subseteq \{e_{ij}\})\Vert\sum_{kl=1\dots
m}\lambda_{kl}(x)e_{kl}\Vert\leq K\},
$$
where $\Vert\sum_{kl=1\dots m}\lambda_{kl}(x)e_{kl}\Vert\leq K$
means $(\forall x_o\in X) \Vert \sum_{kl=1\dots
m}\lambda_{kl}(x_o)e_{kl}\Vert\leq K$. The set $\mathcal{M}$ is a
vector space with point-wise algebraic operations. The map
$\Vert\,\,\, \Vert : \mathcal{M}\to {\bf R}_+$ defined as
$$
\Vert a \Vert = \sup_{\{e_{kl}\}_{kl=1}^n\subseteq
\{e_{ij}\}}\Vert\sum_{kl=1}^n \lambda_{kl}(x)e_{kl}\Vert,
$$
is a norm on the vector space $\mathcal{M}$, where $a\in
\mathcal{M}$ and $a=\sum_{ij\in\Xi}\lambda_{ij}(x)e_{ij}$.

Moreover $\mathcal{M}$ is a von Neumann algebra of
type I$_n$ and $\mathcal{M}=C(X)\otimes M_n({\bf C})$, where the multiplication
is defined as follows $xy=\sum_{ij\in \Xi} \sum_{\xi\in \Xi} \lambda_{i\xi}(x)\mu_{\xi
j}(x)e_{ij}$ \cite{AFN2}.

Let $\mathcal{M}$ be a von Neumann algebra, $\bigtriangleup :\mathcal{M}\to \mathcal{M}$ be a 2-local derivation. Now let
us show that $\bigtriangleup$ is homogeneous. Indeed, for each $x\in \mathcal{M}$, and for
$\lambda \in {\mathbb C}$ there exists a derivation $D_{x,\lambda x}$ such that
$\bigtriangleup(x)=D_{x,\lambda x}(x)$ and $\bigtriangleup(\lambda x)=D_{x,\lambda x}(\lambda x)$. Then
$$
\bigtriangleup(\lambda x)=D_{x,\lambda x}(\lambda x)=\lambda D_{x,\lambda x}(x) =\lambda \bigtriangleup(x).
$$
Hence, $\bigtriangleup$ is homogenous. At the same time,
for each $x\in \mathcal{M}$, there exists a derivation   $D_{x,x^2}$ such that
$\bigtriangleup(x)=D_{x,x^2}(x)$ and $\bigtriangleup(x^2)=D_{x,x^2}(x^2)$.
Then
$$
\bigtriangleup(x^2)=D_{x,x^2}(x^2)=D_{x,x^2}(x)x+xD_{x,x^2}(x) =\bigtriangleup(x)x+x\bigtriangleup(x).
$$

In \cite{Bre} it is proved that any Jordan derivation on a semi-prime algebra is a derivation. Since
$\mathcal{M}$ is semi-prime, the map  $\bigtriangleup$ is a derivation if it is additive.
Therefore, to prove that the 2-local derivation $\bigtriangleup :\mathcal{M}\to \mathcal{M}$ is a derivation
it is sufficient to prove that $\bigtriangleup :\mathcal{M}\to \mathcal{M}$ is additive in the proofs of
theorems 1 and 5.

\section{2-local derivations on von Neumann algebras of type I$_n$
with an infinite cardinal number $n$}

The following theorem is the key  result of this section.

{\bf Theorem 1.} {\it Let $\bigtriangleup :C(X)\otimes M_n({\bf C})\to C(X)\otimes M_n({\bf C})$
be a 2-local derivation. Then $\bigtriangleup$ is a derivation.}

\medskip

First let us prove lemmata which are necessary for the proof of
theorem 1.

Put $\mathcal{M}=C(X)\otimes M_n({\bf C})$, $e_{ij}:={\bf
1}e_{ij}$ for all $i$, $j$, where ${\bf 1}$ is unit of the algebra
$C(X)$. Let $\{a(ij)\}\subset \mathcal{M}$ be the set such that
$$
\bigtriangleup(e_{ij})=a(ij)e_{ij}-e_{ij}a(ij).
$$
for all $i$, $j$,  put $a_{ij}e_{ij}=e_ia(ji)e_j$ for all pairs of
different indexes $i$, $j$ and let
$\{a_{\xi\eta}e_{\xi\eta}\}_{\xi\neq \eta}$ be the set of all such
elements.

\medskip

{\bf Lemma 2.} {\it For any pair $i$, $j$ of different indices the
following equality holds
$$
\bigtriangleup(e_{ij})=\{a_{\xi\eta}e_{\xi\eta}\}_{\xi\neq \eta}e_{ij}-e_{ij}\{a_{\xi\eta}e_{\xi\eta}\}_{\xi\neq \eta}+
a(ij)_{ii}e_{ij}-e_{ij}a(ij)_{jj},    \,\,\,\,\,\,\,(1)
$$
where $a(ij)_{ii}$, $a(ij)_{jj}$ are functions
in $C(X)$ which are the coefficients of the Peirce components $e_{ii}a(ij)e_{ii}$, $e_{jj}a(ij)e_{jj}$.}

{\it Proof.} Let $k$ be an arbitrary index different from $i$, $j$
and let $a(ij,ik)\in \mathcal{M}$ be an element such that
$$
\bigtriangleup(e_{ik})=a(ij,ik)e_{ik}-e_{ik}a(ij,ik) \,\, \text{and}\,\,
\bigtriangleup(e_{ij})=a(ij,ik)e_{ij}-e_{ij}a(ij,ik).
$$
Then
$$
e_{kk}\bigtriangleup(e_{ij})e_{jj}=e_{kk}(a(ij,ik)e_{ij}-e_{ij}a(ij,ik))e_{jj}=
$$
$$
e_{kk}a(ij,ik)e_{ij}-0=e_{kk}a(ik)e_{ij}-e_{kk}e_{ij}\{a_{\xi\eta}e_{\xi\eta}\}_{\xi\neq \eta}e_{jj}=
$$
$$
e_{kk}a_{ki}e_{ij}-e_{kk}e_{ij}\{a_{\xi\eta}e_{\xi\eta}\}_{\xi\neq \eta}e_{jj}=
e_{kk}\{a_{\xi\eta}e_{\xi\eta}\}_{\xi\neq \eta}e_{ij}-e_{kk}e_{ij}\{a_{\xi\eta}e_{\xi\eta}\}_{\xi\neq \eta}e_{jj}=
$$
$$
e_{kk}(\{a_{\xi\eta}e_{\xi\eta}\}_{\xi\neq \eta}e_{ij}-e_{ij}\{a_{\xi\eta}e_{\xi\eta}\}_{\xi\neq \eta})e_{jj}.
$$

Similarly,
$$
e_{kk}\bigtriangleup(e_{ij})e_{ii}=e_{kk}(a(ij,ik)e_{ij}-e_{ij}a(ij,ik))e_{ii}=
$$
$$
e_{kk}a(ij,ik)e_{ij}e_{ii}-0=0-0=
e_{kk}\{a_{\xi\eta}e_{\xi\eta}\}_{\xi\neq \eta}e_{ij}e_{ii}-e_{kk}e_{ij}\{a_{\xi\eta}e_{\xi\eta}\}_{\xi\neq \eta}e_{ii}=
$$
$$
e_{kk}(\{a_{\xi\eta}e_{\xi\eta}\}_{\xi\neq \eta}e_{ij}-e_{ij}\{a_{\xi\eta}e_{\xi\eta}\}_{\xi\neq \eta})e_{ii}.
$$

Let $a(ij,kj)\in \mathcal{M}$ be an element such that
$$
\bigtriangleup(e_{kj})=a(ij,kj)e_{kj}-e_{kj}a(ij,kj)  \,\, \text{and}\,\,
\bigtriangleup(e_{ij})=a(ij,kj)e_{ij}-e_{ij}a(ij,kj).
$$

Then
$$
e_{ii}\bigtriangleup(e_{ij})e_{kk}=e_{ii}(a(ij,kj)e_{ij}-e_{ij}a(ij,kj))e_{kk}=
$$
$$
0-e_{ij}a(ij,kj)e_{kk}=0-e_{ij}a(kj)e_{kk}=0-e_{ij}a_{jk}e_{kk}=
$$
$$
e_{ii}\{a_{\xi\eta}e_{\xi\eta}\}_{\xi\neq \eta}e_{ij}e_{kk}-e_{ij}\{a_{\xi\eta}e_{\xi\eta}\}_{\xi\neq \eta}e_{kk}=
$$
$$
e_{ii}(\{a_{\xi\eta}e_{\xi\eta}\}_{\xi\neq \eta}e_{ij}-e_{ij}\{a_{\xi\eta}e_{\xi\eta}\}_{\xi\neq \eta})e_{kk}.
$$

Also we have
$$
e_{jj}\bigtriangleup(e_{ij})e_{kk}=e_{jj}(a(ij,kj)e_{ij}-e_{ij}a(ij,kj))e_{kk}=
$$
$$
0-0=e_{jj}\{a(ij)\}_{i\neq j}e_{ij}e_{kk}-e_{jj}e_{ij}\{a(ij)\}_{i\neq j}e_{kk}=
$$
$$
e_{jj}(\{a_{\xi\eta}e_{\xi\eta}\}_{\xi\neq \eta}e_{ij}-e_{ij}\{a_{\xi\eta}e_{\xi\eta}\}_{\xi\neq \eta})e_{kk},
$$

$$
e_{ii}\bigtriangleup(e_{ij})e_{ii}=e_{ii}(a(ij)e_{ij}-e_{ij}a(ij))e_{ii}=
$$
$$
0-e_{ij}a(ij)e_{ii}=0-e_{ij}a(ij)e_{ii}=0-e_{ij}a_{ji}e_{ii}=
$$
$$
e_{ii}\{a_{\xi\eta}e_{\xi\eta}\}_{\xi\neq \eta}e_{ij}e_{ii}-e_{ij}\{a_{\xi\eta}e_{\xi\eta}\}_{\xi\neq \eta}e_{ii}=
$$
$$
e_{ii}(\{a_{\xi\eta}e_{\xi\eta}\}_{\xi\neq \eta}e_{ij}-e_{ij}\{a_{\xi\eta}e_{\xi\eta}\}_{\xi\neq \eta})e_{ii}.
$$

$$
e_{jj}\bigtriangleup(e_{ij})e_{jj}=e_{jj}(a(ij)e_{ij}-e_{ij}a(ij))e_{jj}=
$$
$$
e_{jj}a(ij)e_{ij}-0=e_{jj}a_{ji}e_{ij}-0=
$$
$$
e_{jj}\{a_{\xi\eta}e_{\xi\eta}\}_{\xi\neq \eta}e_{ij}-e_{jj}e_{ij}\{a_{\xi\eta}e_{\xi\eta}\}_{\xi\neq \eta}e_{jj}=
$$
$$
e_{jj}(\{a_{\xi\eta}e_{\xi\eta}\}_{\xi\neq \eta}e_{ij}-e_{ij}\{a_{\xi\eta}e_{\xi\eta}\}_{\xi\neq \eta})e_{jj}.
$$

Hence the equality (1) holds. $\triangleright$

We take elements of the sets $\{\{e_{i\xi}\}_\xi\}_i$ and $\{\{e_{\xi j}\}_\xi\}_j$ in
pairs $(\{e_{\alpha\xi}\}_\xi,\{e_{\xi \beta}\}_\xi)$ such that $\alpha\neq \beta$.
Then using the set $\{(\{e_{\alpha\xi}\}_\xi,\{e_{\xi\beta}\}_\xi)\}$ of such pairs
we get the set $\{e_{\alpha\beta}\}$.

Let $x_o=\{e_{\alpha\beta}\}$ be a set $\{v_{ij}e_{ij}\}_{ij}$
such that for all $i$, $j$ if $(\alpha,\beta)\neq (i,j)$ then
$v_{ij}=0$ else $v_{ij}=1$. Then $x_o\in \mathcal{M}$. Fix
different indices $i_o$, $j_o$. Let $c\in \mathcal{M}$ be an
element such that
$$
\bigtriangleup(e_{i_oj_o})=ce_{i_oj_o}-e_{i_oj_o}c \,\,
\text{and}\,\, \bigtriangleup(x_o)=cx_o-x_oc.
$$

Put $c=\{c_{ij}e_{ij}\}\in \mathcal{M}$ and $\bar{a}=\{a_{ij}e_{ij}\}_{i\neq j}\cup \{a_{ii}e_{ii}\}$, where
$\{a_{ii}e_{ii}\}=\{c_{ii}e_{ii}\}$.

\medskip

{\bf Lemma 3.} {\it Let $\xi$, $\eta$ be arbitrary different
indices, and let $b\in \mathcal{M}$ be an element such that
$$
\bigtriangleup(e_{\xi\eta})=be_{\xi\eta}-e_{\xi\eta}b \,\,
\text{and}\,\, \bigtriangleup(x_o)=bx_o-x_ob.
$$
Then $c_{\xi\xi}-c_{\eta\eta}=b_{\xi\xi}-b_{\eta\eta}$.}

{\it Proof.} We have that there exist $\bar{\alpha}$, $\bar{\beta}$ such that
$e_{\xi\bar{\alpha}}$, $e_{\bar{\beta}\eta}\in \{e_{\alpha\beta}\}$ (or $e_{\bar{\alpha}\eta}$,
$e_{\xi\bar{\beta}}\in \{e_{\alpha\beta}\}$, or $e_{\bar{\alpha},\bar{\beta}}\in \{e_{\alpha\beta}\}$), and there
exists a chain of pairs of indexes $(\hat{\alpha}, \hat{\beta})$ in $\Omega$, where
$\Omega=\{(\check{\alpha},\check{\beta}): e_{\check{\alpha},\check{\beta}}\in \{e_{\alpha\beta}\}\}$,
connecting pairs $(\xi, \bar{\alpha})$, $(\bar{\beta},\eta)$ i.e.,
$$
(\xi,\bar{\alpha}), (\bar{\alpha}, \xi_1), (\xi_1, \eta_1),\dots ,(\eta_2, \bar{\beta}), (\bar{\beta},\eta).
$$
Then
$$
c_{\xi\xi}-c_{\bar{\alpha}\bar{\alpha}}=b_{\xi\xi}-b_{\bar{\alpha}\bar{\alpha}},
c_{\bar{\alpha}\bar{\alpha}}-c_{\xi_1\xi_1}=b_{\bar{\alpha}\bar{\alpha}}-b_{\xi_1\xi_1},
$$
$$
c_{\xi_1\xi_1}-c_{\eta_1\eta_1}=b_{\xi_1\xi_1}-b_{\eta_1\eta_1},\dots ,
c_{\eta_2\eta_2}-c_{\bar{\beta}\bar{\beta}}=b_{\eta_2\eta_2}-b_{\bar{\beta}\bar{\beta}},
c_{\bar{\beta}\bar{\beta}}-c_{\eta\eta}=b_{\bar{\beta}\bar{\beta}}-b_{\eta\eta}.
$$
Hence
$$
c_{\xi\xi}-b_{\xi\xi}=c_{\bar{\alpha}\bar{\alpha}}-b_{\bar{\alpha}\bar{\alpha}},
c_{\bar{\alpha}\bar{\alpha}}-b_{\bar{\alpha}\bar{\alpha}}=c_{\xi_1\xi_1}-b_{\xi_1\xi_1},
$$
$$
c_{\xi_1\xi_1}-b_{\xi_1\xi_1}=c_{\eta_1\eta_1}-b_{\eta_1\eta_1},\dots ,
c_{\eta_2\eta_2}-b_{\eta_2\eta_2}=c_{\bar{\beta}\bar{\beta}}-b_{\bar{\beta}\bar{\beta}},
c_{\bar{\beta}\bar{\beta}}-b_{\bar{\beta}\bar{\beta}}=c_{\eta\eta}-b_{\eta\eta}.
$$
and $c_{\xi\xi}-b_{\xi\xi}=c_{\eta\eta}-b_{\eta\eta}$, $c_{\xi\xi}-c_{\eta\eta}=b_{\xi\xi}-b_{\eta\eta}$.

Therefore $c_{\xi\xi}-c_{\eta\eta}=b_{\xi\xi}-b_{\eta\eta}$. $\triangleright$

\medskip

{\bf Lemma 4.} {\it Let $x$ be an element of the algebra $\mathcal{M}$. Then
$$
\bigtriangleup(x)=\bar{a}x-x\bar{a},
$$
where $\bar{a}$ is defined as above.}

{\it Proof.}
Let $d(ij)\in \mathcal{M}$ be an element such that
$$
\bigtriangleup(e_{ij})=d(ij)e_{ij}-e_{ij}d(ij) \,\, \text{and}\,\,
\bigtriangleup(x)=d(ij)x-xd(ij)
$$
and $i\neq j$. Then
$$
\bigtriangleup(e_{ij})=d(ij)e_{ij}-e_{ij}d(ij)=
$$
$$
e_{ii}d(ij)e_{ij}-e_{ij}d(ij)e_{jj}+
(1-e_{ii})d(ij)e_{ij}-e_{ij}d(ij)(1-e_{jj})=
$$
$$
a(ij)_{ii}e_{ij}-e_{ij}a(ij)_{jj}+
\{a_{\xi\eta}e_{\xi\eta}\}_{\xi\neq \eta}e_{ij}-e_{ij}\{a_{\xi\eta}e_{\xi\eta}\}_{\xi\neq \eta}
$$
for all $i$, $j$ by lemma 2.

Since
$$
e_{ii}d(ij)e_{ij}-e_{ij}d(ij)e_{jj}=a(ij)_{ii}e_{ij}-e_{ij}a(ij)_{jj}
$$
we have
$$
(1-e_{ii})d(ij)e_{ii}=\{a_{\xi\eta}e_{\xi\eta}\}_{\xi\neq \eta}e_{ii},
$$
$$
e_{jj}d(ij)(1-e_{jj})=e_{jj}\{a_{\xi\eta}e_{\xi\eta}\}_{\xi\neq \eta}
$$
for all different $i$ and $j$.

Let $b\in \mathcal{M}$ be an element such that
$$
\bigtriangleup(e_{ij})=be_{ij}-e_{ij}b \,\, \text{and}\,\,
\bigtriangleup(x_o)=bx_o-x_ob.
$$
Then $b_{ii}-b_{jj}=c_{ii}-c_{jj}$ by lemma 3. We have
$b_{ii}-b_{jj}=d(ij)_{ii}-d(ij)_{jj}$ since
$$
be_{ij}-e_{ij}b=d(ij)e_{ij}-e_{ij}d(ij).
$$
Hence
$$
c_{ii}-c_{jj}=d(ij)_{ii}-d(ij)_{jj},
c_{jj}-c_{ii}=d(ij)_{jj}-d(ij)_{ii}.
$$

Therefore we have
$$
e_{jj}\bigtriangleup(x)e_{ii}=e_{jj}(d(ij)x-xd(ij))e_{ii}=
$$
$$
e_{jj}d(ij)(1-e_{jj})xe_{ii}+
e_{jj}d(ij)e_{jj}xe_{ii}-e_{jj}x(1-e_{ii})d(ij)e_{ii}-e_{jj}xe_{ii}d(ij)e_{ii}=
$$
$$
e_{jj}\{a_{\xi\eta}e_{\xi\eta}\}_{\xi\neq \eta}xe_{ii}-e_{jj}x\{a_{\xi\eta}e_{\xi\eta}\}_{\xi\neq \eta}e_{ii}+
e_{jj}d(ij)e_{jj}xe_{ii}-e_{jj}xe_{ii}d(ij)e_{ii}=
$$
$$
e_{jj}\{a_{\xi\eta}e_{\xi\eta}\}_{\xi\neq
\eta}xe_{ii}-e_{jj}x\{a_{\xi\eta}e_{\xi\eta}\}_{\xi\neq
\eta}e_{ii}+ c_{jj}e_{jj}xe_{ii}-e_{jj}xe_{ii}c_{ii}e_{ii}=
$$
$$
e_{jj}\{a_{\xi\eta}e_{\xi\eta}\}_{\xi\neq
\eta}xe_{ii}-e_{jj}x\{a_{\xi\eta}e_{\xi\eta}\}_{\xi\neq
\eta}e_{ii}+
$$
$$
e_{jj}(\sum_\xi a_{\xi\xi}e_{\xi\xi})xe_{ii}-e_{jj}x(\sum_\xi
a_{\xi\xi}e_{\xi\xi})e_{ii}=
$$
$$
e_{jj}\{a_{\xi\eta}e_{\xi\eta}\}xe_{ii}-e_{jj}x\{a_{\xi\eta}e_{\xi\eta}\}e_{ii}=
e_{jj}(\bar{a}x-x\bar{a})e_{ii}.
$$

Let $d(ii)$, $v$, $w\in \mathcal{M}$ be elements such that
$$
\bigtriangleup(e_{ii})=d(ii)e_{ii}-e_{ii}d(ii) \,\, \text{and}\,\,
\bigtriangleup(x)=d(ii)x-xd(ii),
$$
$$
\bigtriangleup(e_{ii})=ve_{ii}-e_{ii}v,
\bigtriangleup(e_{ij})=ve_{ij}-e_{ij}v,
$$
and
$$
\bigtriangleup(e_{ii})=we_{ii}-e_{ii}w,
\bigtriangleup(e_{ji})=we_{ji}-e_{ji}w.
$$
Then
$$
(1-e_{ii})a(ij)e_{ii}=(1-e_{ii})ve_{ii}=(1-e_{ii})d(ii)e_{ii},
$$
and
$$
e_{ii}a(ji)(1-e_{ii})=e_{ii}w(1-e_{ii})=e_{ii}d(ii)(1-e_{ii}).
$$
By lemma 2
$$
\bigtriangleup(e_{ij})=a(ij)e_{ij}-e_{ij}a(ij)=
$$
$$
\{a_{\xi\eta}e_{\xi\eta}\}_{\xi\neq
\eta}e_{ij}-e_{ij}\{a_{\xi\eta}e_{\xi\eta}\}_{\xi\neq \eta}+
a(ij)_{ii}e_{ij}-e_{ij}a(ij)_{jj}
$$
and
$$
(1-e_{ii})a(ij)e_{ii}=\{a_{\xi\eta}e_{\xi\eta}\}_{\xi\neq
\eta}e_{ii}.
$$
Similarly
$$
e_{ii}a(ji)(1-e_{ii})=e_{ii}\{a_{\xi\eta}e_{\xi\eta}\}_{\xi\neq
\eta}.
$$

Hence
$$
e_{ii}\bigtriangleup(x)e_{ii}=e_{ii}(d(ii)x-xd(ii))e_{ii}=
$$
$$
e_{ii}d(ii)(1-e_{ii})xe_{ii}+
e_{ii}d(ii)e_{ii}xe_{ii}-e_{ii}x(1-e_{ii})d(ii)e_{ii}-e_{ii}xe_{ii}d(ii)e_{ii}=
$$
$$
e_{ii}a(ji)(1-e_{ii})xe_{ii}+
e_{ii}d(ii)e_{ii}xe_{ii}-e_{ii}x(1-e_{ii})a(ij)e_{ii}-e_{ii}xe_{ii}d(ii)e_{ii}=
$$
$$
e_{ii}\{a_{\xi\eta}e_{\xi\eta}\}_{\xi\neq \eta}xe_{ii}-e_{ii}x\{a_{\xi\eta}e_{\xi\eta}\}_{\xi\neq \eta}e_{ii}+
e_{ii}d(ii)e_{ii}xe_{ii}-e_{ii}xe_{ii}d(ii)e_{ii}=
$$
$$
e_{ii}\{a_{\xi\eta}e_{\xi\eta}\}_{\xi\neq
\eta}xe_{ii}-e_{ii}x\{a_{\xi\eta}e_{\xi\eta}\}_{\xi\neq
\eta}e_{ii}+c_{ii}e_{ii}xe_{ii}-e_{ii}xc_{ii}e_{ii}=
$$
$$
e_{ii}\{a_{\xi\eta}e_{\xi\eta}\}_{\xi\neq
\eta}xe_{ii}-e_{ii}x\{a_{\xi\eta}e_{\xi\eta}\}_{\xi\neq
\eta}e_{ii}+
$$
$$
e_{ii}(\sum_\xi a_{\xi\xi}e_{\xi\xi})xe_{ii}-e_{ii}x(\sum_\xi
a_{\xi\xi}e_{\xi\xi})e_{ii}=
$$
$$
e_{ii}\{a_{\xi\eta}e_{\xi\eta}\}xe_{ii}-e_{ii}x\{a_{\xi\eta}e_{\xi\eta}\}e_{ii}=e_{ii}(\bar{a}x-x\bar{a})e_{ii}.
$$

Hence
$$
\bigtriangleup(x)=\bar{a}x-x\bar{a}
$$
for all $x\in \mathcal{M}$.
$\triangleright$

\medskip

{\it Proof of theorem 1.} Let $V=\{\{\lambda_{ij}e_{ij}\}_{ij}: \{\lambda_{ij}\}\subset C(X)\}$ (the set of all
infinite $n\times n$-dimensional function-valued matrices). Then $V$ is a vector space with componentwise
algebraic operations and $\mathcal{M}$ is a vector subspace of $V$.

By lemma 4 $\bigtriangleup(e_{ii})=\bar{a}e_{ii}-e_{ii}\bar{a}\in \mathcal{M}$. Hence
$$
\sum_\xi a_{\xi i}e_{\xi i}-\sum_\xi a_{i\xi}e_{i\xi}\in \mathcal{M}.
$$
Then
$$
e_{ii}(\sum_\xi a_{\xi i}e_{\xi i}-\sum_\xi a_{i\xi}e_{i\xi})=
a_{ii}e_{ii}-\sum_\xi a_{i\xi}e_{i\xi}\in \mathcal{M}
$$
and
$$
(\sum_\xi a_{\xi i}e_{\xi i}-\sum_\xi a_{i\xi}e_{i\xi})e_{ii}=
\sum_\xi a_{\xi i}e_{\xi i}-a_{ii}e_{ii}\in \mathcal{M}.
$$
Therefore $\sum_\xi a_{\xi i}e_{\xi i}$, $\sum_\xi a_{i\xi}e_{i\xi}\in \mathcal{M}$ i.e.,
$\bar{a}e_{ii}, e_{ii}\bar{a}\in \mathcal{M}$. Hence $e_{ii}\bar{a}x, x\bar{a}e_{ii}\in \mathcal{M}$ for any $i$
and
$$
\bar{a}x, x\bar{a}\in V
$$
for any element $x=\{x_{ij}e_{ij}\}\in \mathcal{M}$, i.e.,
$$
\sum_\xi a_{i\xi}x_{\xi j}e_{ij}, \sum_\xi x_{i\xi}a_{\xi j}e_{ij}\in {\mathbb C}e_{ij}
$$
for all $i$, $j$. Therefore for all $x$, $y\in \mathcal{M}$ we
have that the elements $\bar{a}x$, $x\bar{a}$, $\bar{a}y$,
$y\bar{a}$, $\bar{a}(x+y)$, $(x+y)\bar{a}$ belong to $V$. Hence
$$
\bigtriangleup(x+y)=\bigtriangleup(x)+\bigtriangleup(y)
$$
by lemma 4.

Similarly for all $x$, $y\in \mathcal{M}$ we have
$$
(\bar{a}x+x\bar{a})y=\bar{a}xy-x\bar{a}y\in \mathcal{M}, \bar{a}xy=\bar{a}(xy)\in V.
$$
Then $x\bar{a}y=\bar{a}xy-(\bar{a}x-x\bar{a})y$ and $x\bar{a}y\in V$.
Therefore
$$
\bar{a}(xy)-(xy)\bar{a}=\bar{a}xy-x\bar{a}y+x\bar{a}y-xy\bar{a}=(\bar{a}x-x\bar{a})y+x(\bar{a}y-y\bar{a}).
$$
Hence
$$
\bigtriangleup(xy)=\bigtriangleup(x)y+x\bigtriangleup(y)
$$
by lemma 4.
Now we show that $\bigtriangleup$ is homogeneous. Indeed, for each $x\in \mathcal{M}$, and for
$\lambda \in {\mathbb C}$ there exists a derivation $D_{x,\lambda x}$ such that
$\bigtriangleup(x)=D_{x,\lambda x}(x)$ and $\bigtriangleup(\lambda x)=D_{x,\lambda x}(\lambda x)$. Then
$$
\bigtriangleup(\lambda x)=D_{x,\lambda x}(\lambda x)=\lambda D_{x,\lambda x}(x) =\lambda \bigtriangleup(x).
$$
Hence, $\bigtriangleup$ is homogenous and therefore it is a linear operator and a derivation.
The proof is complete.

$\triangleright$

\section{The main theorem}

{\bf Theorem 5.} {\it Let $M$ be a von Neumann algebra of type $I$
and let $\bigtriangleup :M\to M$ be a 2-local derivation. Then
$\bigtriangleup$ is a derivation.}

{\it Proof.} We have that
$$
M=\sum_j^\oplus M_{I_{n_j}},
$$
where $M_{I_{n_j}}$ is a von Neumann algebra of type $I_{n_j}$, $n_j$ is a cardinal number for any $j$.
Let $x_j\in M_{I_{n_j}}$ for any $j$ and $x=\sum_j x_j$. Note that
$\bigtriangleup(x_j)\in M_{I_{n_j}}$ for all $x_j\in M_{I_{n_j}}$. Hence
$$
\bigtriangleup\vert_{M_{I_{n_j}}}: M_{I_{n_j}}\to M_{I_{n_j}}
$$
and $\bigtriangleup$ is a 2-local derivation on $M_{I_{n_j}}$. There exists a hyperstonean
compact $X$ such that $M_{I_{n_j}}\cong C(X)\otimes M_{n_j}({\bf C})$. Hence by theorem 1
$\bigtriangleup$ is a derivation on $M_{I_{n_j}}$.

Let $x$ be an arbitrary element of $M$.
Then there exists $d(j)\in M$ such that $\bigtriangleup(x)=d(j)x-xd(j)$,
$\bigtriangleup(x_j)=d(j)x_j-x_jd(j)$ and
$$
z_j\bigtriangleup(x)=z_j(d(j)x-xd(j))=z_j\sum_i (d(j)x_i-x_id(j))=
$$
$$
d(j)x_j-x_jd(j)=\bigtriangleup(x_j),
$$
for all $j$, where $z_j$ is unit of $M_{I_{n_j}}$. Hence
$$
\bigtriangleup(x)=\sum_j z_j\bigtriangleup(x)=\sum_j \bigtriangleup(x_j).
$$
Since $x$ was chosen arbitrarily $\bigtriangleup$ is a derivation on $M$ by the last
equality.

Indeed, let $x,y\in M$. Then
$$
\bigtriangleup(x)+\bigtriangleup(y)=\sum_j \bigtriangleup(x_j)+\sum_j \bigtriangleup(y_j)=
\sum_j [\bigtriangleup(x_j)+\bigtriangleup(y_j)]=
$$
$$
\sum_j \bigtriangleup(x_j+y_j)=\sum_j z_j\bigtriangleup(x+y)=
\bigtriangleup(x+y).
$$

Similarly,
$$
\bigtriangleup(xy)=\sum_j \bigtriangleup(x_jy_j)=\sum_j [\bigtriangleup(x_j)y_j+x_j\bigtriangleup(y_j)]=
$$
$$
\sum_j \bigtriangleup(x_j)y_j+\sum_j x_j\bigtriangleup(y_j)=
\sum_j \bigtriangleup(x_j) \sum_j y_j+\sum_j x_j \sum_j \bigtriangleup(y_j)=
$$
$$
\bigtriangleup(x)y+x\bigtriangleup(y).
$$
By the proof of the previous theorem $\bigtriangleup$ is homogenous.
Hence $\bigtriangleup$ is a linear operator and a derivation.
The proof is complete.
$\triangleright$

\end{document}